\documentclass{IOS-Book-Article}

\usepackage{mathptmx}
\usepackage{dsfont}
\usepackage{soul}\setuldepth{article}
\usepackage{times}
\usepackage{soul}
\usepackage{url}
\usepackage[utf8]{inputenc}
\usepackage[small]{caption}
\usepackage{graphicx}
\usepackage{eqnarray,amsmath}
\usepackage{amsthm}
\usepackage{booktabs}
\usepackage{algorithm}
\usepackage{algorithmic}
\usepackage{cite}

\usepackage{bm}
\usepackage{mathrsfs}
\usepackage{amssymb}

\def\hb{\hbox to 11.5 cm{}}

\begin{document}

\pagestyle{headings}
\def\thepage{}

\begin{frontmatter}              % The preamble begins here.

%\pretitle{Pretitle}
\title{An optimization model with\\stochastic variables for flexible production logistics planning}

\markboth{}{Nov 2021\hb}
%\subtitle{Subtitle}

\author[A]{\fnms{Yongkuk} \snm{Jeong}
\thanks{Corresponding Author: Yongkuk Jeong; E-mail: yongkuk@kth.se}},
\author[B]{\fnms{Gianpiero} \snm{Canessa}}
\author[A]{\fnms{Erik} \snm{Flores-García}}
\author[A]{\fnms{Tarun} \snm{Kumar Agrawal}}
and
\author[A]{\fnms{Magnus} \snm{Wiktorsson}}

\runningauthor{B.P. Manager et al.}
\address[A]{Department of Sustainable Production Development, KTH Royal Institute of Technology, Sweden}
\address[B]{Department of Mathematics, KTH Royal Institute of Technology, Sweden}

\begin{abstract}
Production logistics has an important role as a chain that connects the components of the production system. The most important goal of production logistics plans is to keep the flow of the production system well. However, compared to the production system, the level of planning, management, and digitalization of the production logistics system is not high enough, so it is difficult to respond flexibly when unexpected situations occur in the production logistics system. Optimization and heuristic algorithms have been proposed to solve this problem, but due to their inflexible nature, they can only achieve the desired solution in a limited environment. In this paper, the relationship between the production and production logistics system is analyzed and stochastic variables are introduced by modifying the pickup and delivery problem with time windows (PDPTW) optimization model to establish a flexible production logistics plan. This model, taking into account stochastic variables, gives the scheduler a new perspective, allowing them to have new insights based on the mathematical model. However, since the optimization model is still insufficient to respond to the dynamic environment, future research will cover how to derive meaningful results even in a dynamic environment such as a machine learning model.
\end{abstract}

\end{frontmatter}
\markboth{Nov 2021\hb}{Nov 2021\hb}

\section{Introduction}

\begin{figure}[h]
\centering
\includegraphics[width=0.75\textwidth]{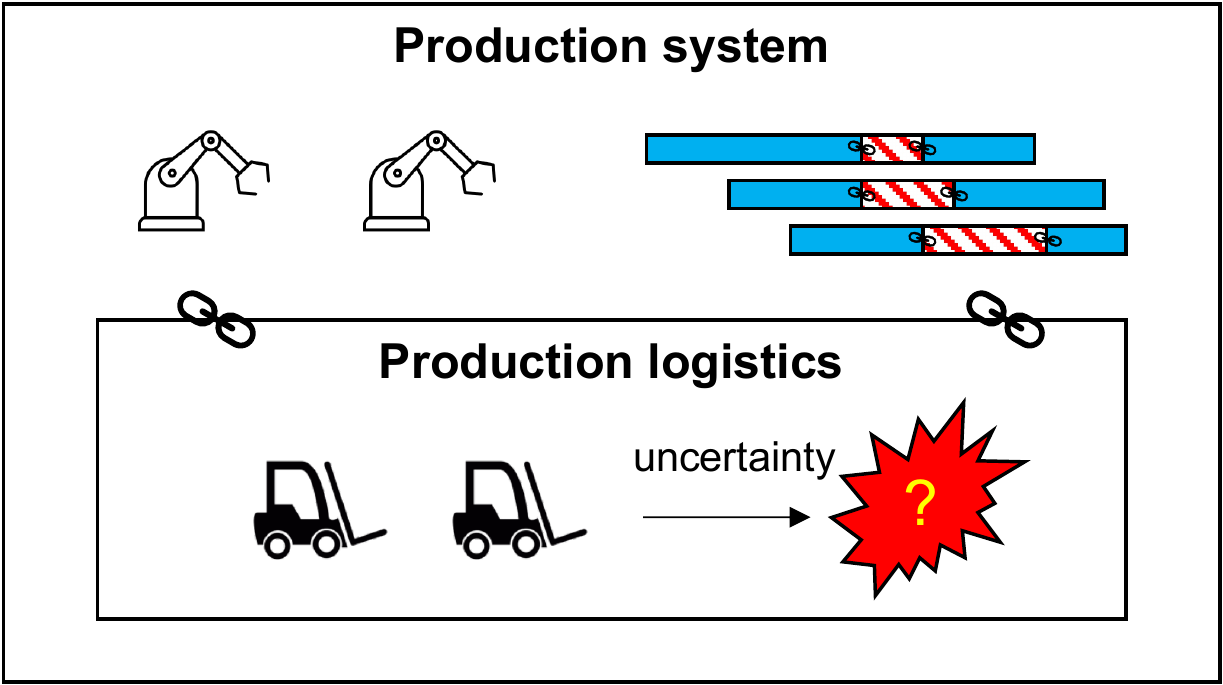}
\caption{Production system and linked production logistics}
\label{fig1:production}
\end{figure}

Production logistics refers to the flow of material and information between components of a production system. The production logistics plan is mostly dependent on the production plan because production logistics is a part that supports the production system. The complexity of the production system is gradually increasing, accordingly, a more sophisticated and flexible plan is required (Figure \ref{fig1:production}). In particular, as the production logistics plan is restricted by compliance with the production plan, a flexible and responsive planning model is required to cope with the situation of the production system. Traditionally, mathematical and heuristic models have been used for production and logistics planning. They have the advantage of being able to find an exact solution or a solution close to it for a given problem but also have the disadvantage of low flexibility for similar problems or new situations.

\begin{figure}[h]
\centering
\includegraphics[width=0.75\textwidth]{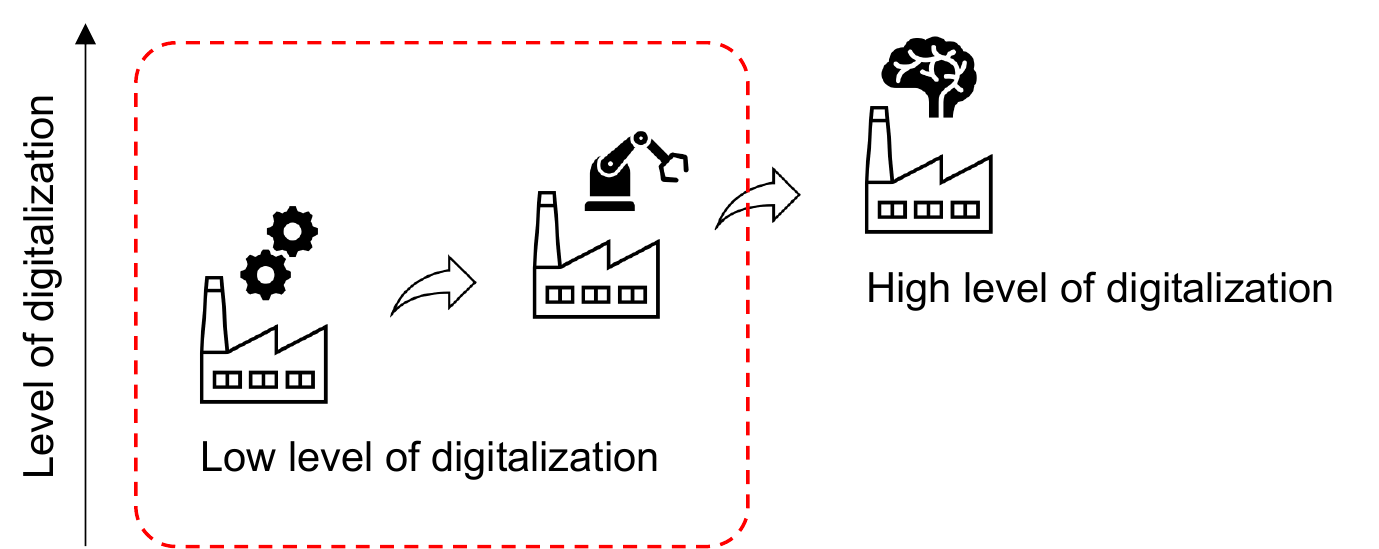}
\caption{Level of digitalization and our interest area}
\label{fig2:target}
\end{figure}

There have been many developments in the field of production system management due to the improvement of the level of digitalization such as automation and Internet-of-Things (IoT). However, since the essential characteristic of production logistics as “production system support”, the area of planning and management of production logistics was not developed as much as the level of planning and management of production systems. This phenomenon can be seen more prominently because the priority of production logistics is lower in an environment with a low level of digitalization (Figure \ref{fig2:target}). As mentioned above, since production logistics plays an important role as a chain that connects the elements of the production system, it is difficult to expect the overall development of the production system without the development of production logistics planning and management. Therefore, we introduce an optimization model that can be used even in a production system environment with a low level of digitalization in this paper. In addition, It is also not only considering the improvement of efficiency of system but also sustainability.

\section{Related works}
\subsection{Production logistics and planning}

Significant research in past has focused on production logistics, production planning and understanding the interlink/interplay between them. In past researchers have proposed several techniques to optimize both the plans to achieve maximum flexibility and efficiency. Particularly, introduction of  Industry 4.0 and smart manufacturing, supported by digital technologies and data analysis methods, has bought a paradigm shift in the area of production logistics planning. Mourtzis suggested a cloud-based cyber-physical system to facilitate adaptive production scheduling with condition-based  maintenance through a cost effective and reliable data acquisition, processing,  and analysis \cite{mourtzis2020simulation}. Alexopoulos demonstrated an Industrial Internet of Things (IIoT) information system, which could offer decision-making support for moving or static operators/supervisors \cite{alexopoulos2018industrial}.

Currently, it is important that the production planning and logistics operations are able to realize the rapid response and communication interaction of high volume of real-time heterogeneous data for large-scale customized product manufacturing \cite{huang2017collaborative}. Production resources including those of production logistics should be able to autonomously and correctly make collective production decisions by forming a comprehensive knowledge system using cloud manufacturing equipment and data analysis to meet customized demands of customers.
Approaches to production logistics path planning have been discussed in various operational scenarios. For instance, \cite{gilboa2006distributed} explained a distributed algorithm through a scenario where initially the layout model of the environment is unknown and the agents physically move to explore new coordinates and vertices.

On the other hand \cite{bhattacharya2011distributed} demonstrated  a distributed optimization algorithm for multi-robot path planning. This algorithm was designed to find paths in a two-dimensional continuous environment based on rendezvous constraints.  In contrast, \cite{chouhan2015dmapp} propose DMAPP - a token-passing algorithm in which each individual agent plans a path through a FastForward planner and then each agent communicates and share its plan with the others based on a priority order to resolve conflicts among the plans.

\subsection{Sustainability in production logistics}

The literature increasingly highlights the need for enhancing sustainability in production logistics including three aspects. First, economic sustainability referring to economic growth without affecting negatively social and environmental aspects. Second, environmental sustainability comprises the responsible interaction with the environment avoiding the degradation of natural resources. Third, social sustainability consists of establishing and managing positive and negative aspects business on staff. Current understanding about sustainability in production logistics points to the need for understanding the trade-off between economic, environmental, and social aspects. However, existing studies in production logistics give precedence to the economic and environmental side, while social sustainability remains understudied.

A number of recent publications underscore advances towards consideration of economic and environmental sustainability in the scheduling and routing in production logistics. For example, \cite{banyai2019smart} describe the in-plant supply process of matrix production and the optimization potential for scheduling material handling. The study applies a sequential black hole–floral pollination heuristic algorithm for minimizing period deviance, energy consumption, route length, and emissions. \cite{liao2019novel} establish an intelligent production scheduling and logistics delivery model with IoT technology to promote green and sustainable development of manufacturing. The study proposes a two-layer optimization applying a genetic algorithm (GA) for minimizing operating costs and carbon emissions when scheduling a production and delivery process. \cite{li2018new} address the problem of automated guided vehicle (AGV) scheduling for transferring materials to CNC machines. The study applies a harmony search (HS) algorithm meta-heuristic for minimizing the standard deviation of the waiting time of CNC material buffers and the total travel distance of the AGV. \cite{yao2020improving} focus on manufacturing disruptions affecting both AGV and machine scheduling. In particular, the study applies a nonlinear mixed integer programming GA to find near-optimal production schedules prioritizing the just-in-time material delivery performance and energy efficiency of the material transportation.

Despite these advances, consideration of social sustainability involving scheduling and routing in production logistics remains scarce. Current efforts either neglect tasks in production logistics or do not take into account scheduling and routing. For example, \cite{otto2017reducing} provide an overview of the existing optimization approaches to assembly line balancing and job rotation scheduling that consider physical ergonomic risks. \cite{lodree2009taxonomy} present a taxonomy for integrating scheduling theory and human factors. \cite{bechtsis2017sustainable} provide a general framework describing the influence of AGVs on economic, environmental, and social sustainability. \cite{halawa2020introduction} analyze tracking of material handling resources for enhancing warehouse safety and operational efficiency.

\section{Production system and production logistics planning}

\subsection{Production logistics planning model}

\begin{figure}[h]
\centering
\includegraphics[width=1.0\textwidth]{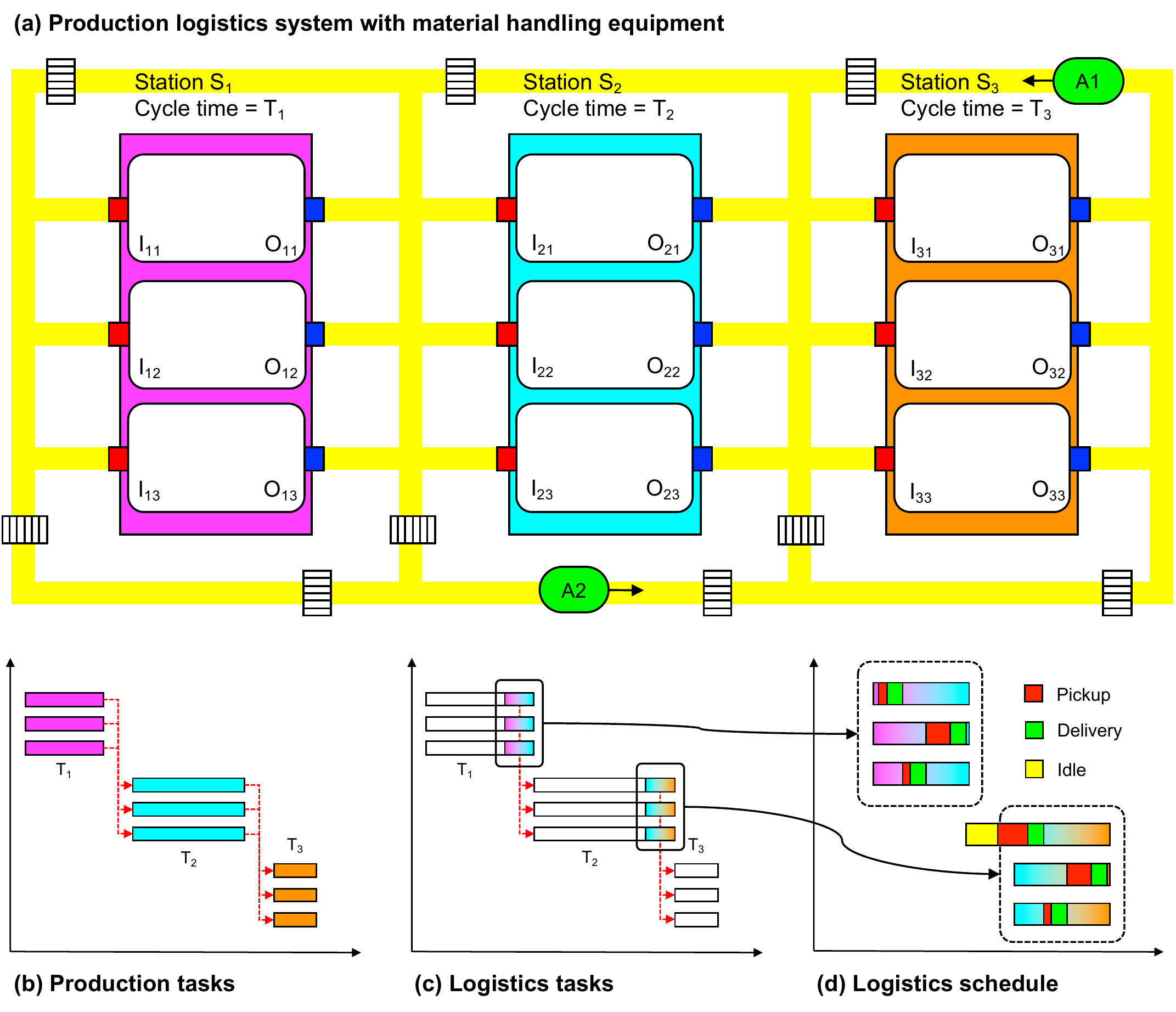}
\caption{Production logistics system and logistics scheduling}
\label{fig3:PL system and logsitics}
\end{figure}

Production system includes production station/unit, material handling equipment (MHE) and as shown in Figure \ref{fig3:PL system and logsitics} (a). Most likely, decision support system for production management focuses on the planning and scheduling of production stations and units. However, the most important goal of production logistics planning is to deliver the right product without interrupting the flow of the production system which is a significant role to achieve the planned throughput of production system.

Therefore MHE will follow the best decision in a given situation rather than follow a planned logistics schedule. There are several reasons for this. First, logistics tasks themselves act as buffers between production tasks. Second, logistics tasks are more likely to experience unexpected disturbances than production tasks. So, if the plan should not be too tight. Third, the production logistics environment is not controlled as much as production stations and units. Therefore, there is a high possibility of encountering an unexpected situation.

Figure \ref{fig3:PL system and logsitics} (b), (c), (d) shows the logistics schedule assigned to production tasks, logistics tasks, and MHE, respectively. Logistics tasks are defined between production tasks (Figure \ref{fig3:PL system and logsitics} (b) and (c)) and they are assigned to the MHE (Figure \ref{fig3:PL system and logsitics} (d)) based on the knowledge and intelligence of the logistics scheduler. Even when there is an urgent task that needs to be handled or an unexpected situation occurs in the production logistics system, the logistics scheduling process will be updated based on the intuition and experience of the logistics scheduler.

The logistics schedule from this process does not guarantee the optimal schedule, and it is not a fundamental solution because the solution is not derived by analyzing the essence of the problem. To handle this problem,  we propose a flexible production logistics planning model by modeling predictable disturbances and situations as stochastic variables. The following section explains this in more detail with a simple example.

\subsection{Sample case}

\begin{figure}[h]
\centering
\includegraphics[width=1.0\textwidth]{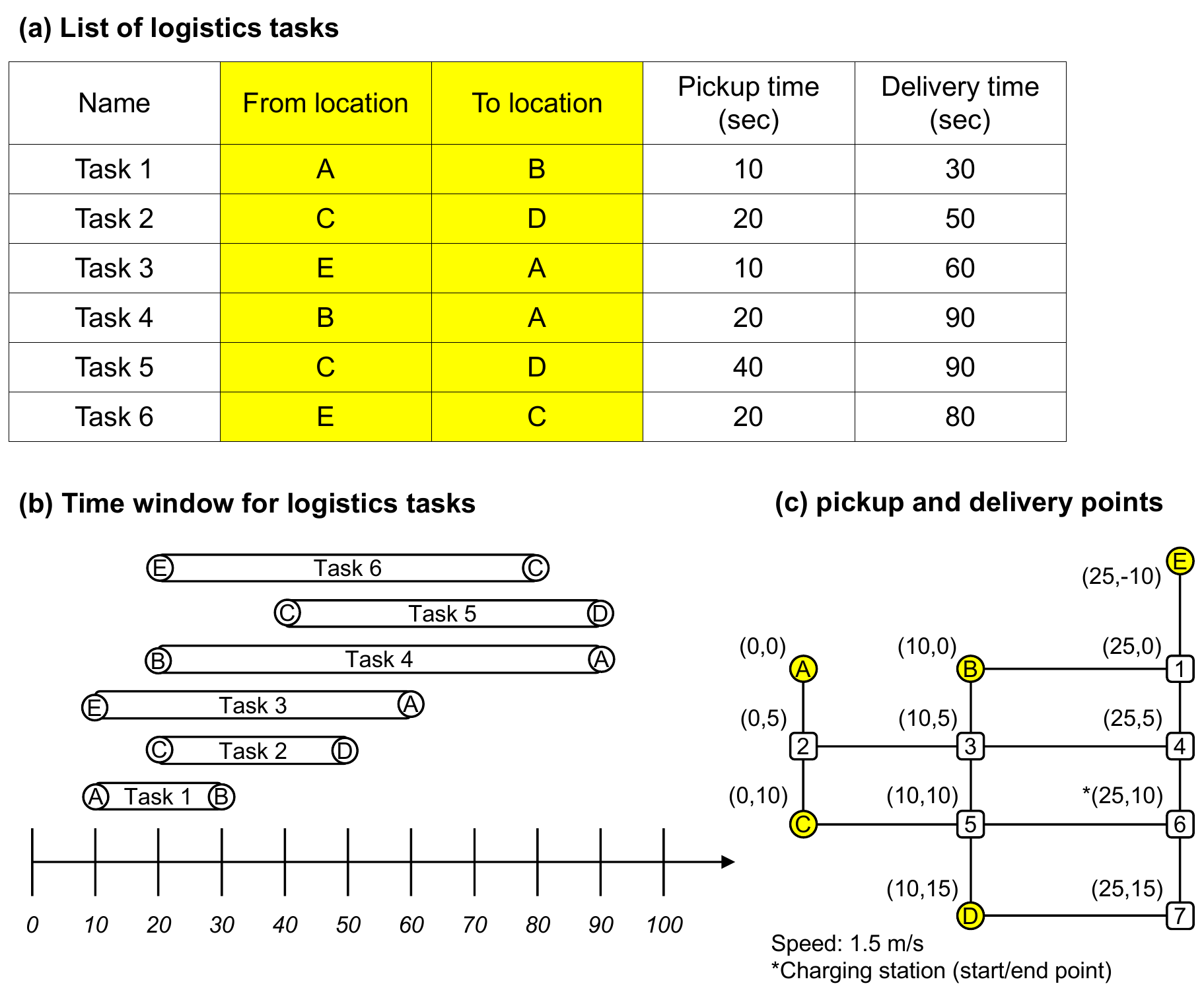}
\caption{Simplified case for production logistics planning}
\label{fig4:example}
\end{figure}

Logistics tasks between production stations can be described with the information as shown in Figure \ref{fig4:example} (a). First, from location indicates the location where the product made and should be picked up. To location means the location where it should be delivered. For example, Task 1 in Figure \ref{fig4:example} (a) will be picked up at Point A and delivered to Point B. Pickup time and delivery time refer to the earliest time that the product will be ready after the production process is completed, and the latest time to deliver the product so that the next production process is not delayed, respectively. Therefore, Task 1 can be picked up after 10 second, and must be delivered to B before 30 second at the latest. It can be presented as a Gantt chart as in Figure \ref{fig4:example} (b).

However, the traveling time between the production stations would be determined differently depending on layout, specification of MHEs. For example, in the layout configured as in Figure \ref{fig4:example} (c), it takes the longest physically to reach from Point A to Point E then other points. In this situation, “A-2-3-4--1-E” or “A-2-3-5-6-4-1-E” could be selected as a delivery route. As a summary, the logistics task can provide basic information, but the time taken to perform the actual logistics tasks can be derived differently depending on the situation. In the next chapter, we propose an optimization model that considers stochastic variables to consider these uncertainties in production logistics environment

\section{An optimization model with stochastic variables}
Our proposal consists of a modified version of the pickup and delivery problem with time windows (PDPTW) model from \cite{battarra2014chapter}. The solution obtained is the order in which each order must be picked up and delivered by each vehicle capable of transporting goods.

\subsection{Modified PDPTW}
First, we will present the general model with minor modifications in its deterministic version.

Let $G(V,E)$ is a complete graph, where $V = \{0, 1, \dots, 2n+1\}$ and $E$ are the node and arc set respectively, which represents our context where each node will be pickup and/or delivery points, with a total of $n$ nodes. Let $P = \{1, \dots, n\}$ and $D = \{1+n, \dots, 2n\}$ be the pickup and delivery node sets respectively. A time window $[a_i, b_i]$ is imposed for each node $i \in N$, where $a_i$ and $b_i$ represent the earliest and latest time respectively at which the vehicle needs to start the pickup service from node $i$. Let $d_e$ be the travel time associated to each arc $e \in E$
%Assume the $c_e$ includes all costs associated to sustainability and finances (GC: No idea how to phrase the fact we may need money here, monetary?)
where it also considers service times at the start node $i \in V$ and target node $j \in V$ of arc $e = (i,j) \in E$.

The PDPTW also takes weight and volume capacities into account, however in our problem there is no need to take these factors into account: each vehicle has enough capacity to serve every task, but can only work on one task a time.

Let $x_{k,e} \in \mathds{B}, \forall k \in K, \forall e \in E$ be a binary decision variable (i.e. $\mathds {B} := \{0,1\}$) for every vehicle $k$ to take a route with arc $e$. Let $w_{k,i} \ge 0, \forall k \in K, \forall i \in V$ be a decision variable representing the time when vehicle $k$ reaches node $i$. Finally, the model is as follows:

\begin{eqnarray}
    DM(d) = & \min_{x,w} & f(x) \label{form:FS_OF} \\
    & \text{s.t.} & \sum_{k \in K} \sum_{j \in V} x_{k,(i,j)} = 1, \forall i \in P, \label{form:FS_c1} \\
    & & \sum_{j \in V} x_{k, (i,j)} - \sum_{j \in V: (j,t) \in E} x_{k,(1+n,j)} = 0, \forall i \in P, \forall k \in K, \label{form:FS_c2}\\
    & & \sum_{j \in V} x_{k,(0,j)} = 1, \forall k \in K, \label{form:FS_c3}\\
    & & \sum_{j \in V} x_{k,(i,2n+1)} = 1, \forall k \in K, \label{form:FS_c4}\\
    & & \sum_{j \in V} x_{k, (i,j)} - \sum_{j \in V: (j,t) \in E} x_{k,(j,i)} = 0, \forall i \in P \cup D, \forall k \in K, \label{form:FS_c5}\\
    & & w_{k,j} \ge w_{k, i} + d_{(i,j)} - M (1 - x_{k,(i,j)}), \forall i \in N, \forall j \in V, \forall k \in K, \label{form:FS_c6}\\
    & & w_{k,i} + d_{(i,i+n)} \le w_{k, i+n}, \forall i \in P, \forall k \in K, \label{form:FS_c7}\\
    & & a_i \le w_{k,i} \le b_i, \forall i \in V, \forall k \in K, \label{form:FS_c8}\\
    & & x_{k,e} \in \mathbb{B}, \forall e \in E, \forall k \in K, \label{form:FS_c10}
\end{eqnarray}

The multi-objective function \eqref{form:FS_OF} minimizes a convex function $f_1:\mathds{R}^{|K| \times |E|} \rightarrow \mathds{R}$. Constraints \eqref{form:FS_c1} and \eqref{form:FS_c2} guarantees that each pair pickup/destiny nodes are visited by the same vehicle. Constraints \eqref{form:FS_c3} and \eqref{form:FS_c4} make sure that every vehicle starts at node 0 and finishes at node $2n+1$, which have been selected as source and terminal. Please note these last two constraints are optional, however if used and there are vehicles not used in any particular route, they will still travel from source to terminal. Constraint \eqref{form:FS_c5} ensures that the vehicle which enters node $i$ must leave the node. The order and time of visitation of every node is given by constraint \eqref{form:FS_c6}, where $M$ is a large enough number. Constraint \eqref{form:FS_c7} guarantees that the pickup node is visited before the delivery node, while \eqref{form:FS_c8} imposes the time windows for each node. Finally, constraint \eqref{form:FS_c10} set the feasible space for $x$.

\subsection{Stochastic modified PDPTW}
A new proposed modification is the nature of the each vector $d_e$, where instead of a deterministic value, let $\tilde{d}_e \sim D_e(\lambda_e)$ be a random variable where $D_e(\cdot)$ is a probability distribution and $\lambda_e$ is the parameter vector with appropriate dimension. For our case, let $\Omega$ be the set of possible scenarios for $\tilde{d}_e, \forall e \in E$, where the probability of each scenario $\omega \in \Omega$ is $p_\omega$. Therefore, constraints \eqref{form:FS_c6} and \eqref{form:FS_c7} need to be modified to cope with the random variable:

\begin{eqnarray*}
    w_{k,j} \ge w_{k, i} + \tilde{d}_{(i,j)} - M (1 - x_{k,(i,j)}), \forall i \in N, \forall j \in V, \forall k \in K, \label{form:FS_c6m}\\
    w_{k,i} + \tilde{d}_{(i,i+n)} \le w_{k, i+n}, \forall i \in P, \forall k \in K, \label{form:FS_c7m}
\end{eqnarray*}

Since these constraint can not work using random variables in this state, we will expand the definition of the variables proposed: Each route described by $x$ will be fixed across all scenarios, while we will give freedom to the model to modify its pickup/delivery time $w$, which will depend on the scenario. Then, let $w_{k,i,\omega}, \forall k \in K, \forall i \in V, \forall \omega \in \Omega$ be the node service times for the scenario $\omega$. Now, the extended model constraint will be formulated using a \emph{chance-constraint} to bound the probability of the time windows compliance as follows:

\begin{equation} \label{}
    \mathds{P}\left[\begin{array}{l}
        w_{k,j} \ge w_{k, i} + \tilde{d}_{(i,j)} - M (1 - x_{k,(i,j)}), \forall i \in N, \forall j \in V, \forall k \in K,\\
        w_{k,i} + \tilde{d}_{(i,i+n)} \le w_{k, i+n}, \forall i \in P, \forall k \in K,\\
        a_i \le w_{k,i}, \quad w_{k,i} \le b_i, \forall i \in V, \forall k \in K.
    \end{array}
    \right] \ge 1 - \alpha,
\end{equation}

where $0 \le alpha < 1$ is defined as the reliability of the system: If $\alpha = 0$, then we assume robustness and we may solve the problem with the worst case scenario possible, and any value $\alpha > 0$ will force the solution to be robust for at least that proportion of scenarios provided. Unfortunately, chance-constraints do not offer an intuitive and out-of-the-box implementation, therefore a reformulation is needed in order to obtain a proper model to solve. We will assume $|\Omega| < \infty$, then:

\begin{eqnarray}
    SPM(d) \label{form:SPM} = & \min_{x,w,c} & f(x)\label{form:SPM_OF} \\
    & \text{s.t.} & \sum_{k \in K} \sum_{j \in V} x_{k,(i,j)} = 1, \forall i \in P, \label{form:SPM_c1}\\
    & & \sum_{j \in V} x_{k, (i,j)} - \sum_{j \in V: (j,t) \in E} x_{k,(1+n,j)} = 0, \forall i \in P, \forall k \in K, \label{form:SPM_c2}\\
    & & \sum_{j \in V} x_{k,(0,j)} = 1, \forall k \in K, \label{form:SPM_c3}\\
    & & \sum_{j \in V} x_{k,(i,2n+1)} = 1, \forall k \in K, \label{form:SPM_c4}\\
    & & \sum_{j \in V} x_{k, (i,j)} - \sum_{j \in V: (j,t) \in E} x_{k,(j,i)} = 0, \forall i \in P \cup D, \forall k \in K, \label{form:SPM_c5}\\
    & & w_{k,j,\omega} \ge w_{k, i, \omega} + d_{(i,j), \omega} - M_1 (1 - x_{k,(i,j)}) - M_2 z_\omega, \notag  \\
    & & \forall i \in N, \forall j \in V, \forall k \in K, \forall \omega \in \Omega,\\
    & & w_{k,i,\omega} + d_{(i,i+n), \omega} \le w_{k, i+n, \omega} + M_3 z_\omega, \notag \\
    & & \forall i \in P, \forall k \in K, \forall \omega \in \Omega,\\
    & & a_i z_\omega \le w_{k,i,\omega} \le b_i + M_4 z_\omega, \forall i \in V, \forall k \in K, \omega \in \Omega,\\
    & & \sum_{\forall \omega \in \Omega} p_\omega z_\omega \le \alpha, \label{form:SPM_c10}\\
    & & x_{k,e} \in \mathbb{B}, \forall e \in E, \forall k \in K, \label{form:SPM_c11}\\
    & & z_\omega \in \mathds{B}, \forall \omega \in \Omega \label{form:SPM_c13},
\end{eqnarray}

where $M_1,\dots, M_5$ are sufficiently large numbers, $z_\omega, \forall \omega \in \Omega$ is a set of binary variables which \emph{select} which scenarios will be \emph{ignored} from the set $\Omega$. Constraint \eqref{form:SPM_c10} bounds how many scenarios can be ignored by bounding the weighted sum of the probability of each scenario ignored by $\alpha$.

If the assumption over the size of $\Omega$ is relaxed, then a sample average approximation (SAA) approach can be used as proposed in \cite{kleywegt2002sample,luedtke2008sample, pagnoncelli2009sample}: let $\hat{\Omega} \subset \Omega, |\hat{\Omega}| < \infty$ and $p_\omega = \frac{1}{\hat{\Omega}}, \forall \omega \in \hat{\Omega}$. If the samples size in $\hat{\Omega}$ is sufficient, then we can bound the probability of model \eqref{form:SPM} to converge to the true optimal solution.

\section{Experimental results}
We tested our model using the example in Figure§ \ref{fig4:example}, while adding a random distribution to travel times. In our example, in order to obtain the distances, we used the graph already presented and then assumed a speed of 1.5 m/s. Let $c_e, \forall e \in E$ be the traveling distance for every arc $e$ and $v = 1.5$, then the deterministic travel time is calculated as follows: $d_e = c_e / v, \forall e \in E$.

A truncated normal distribution with mean 1 and variance 0.25 was utilized as multiplier to the deterministic travel distance in order to create uncertainty: $\tilde{d}_e \sim \max\{N(1, 0.25), 0\}d_e$. With these parameter, the probability of obtaining a negative speed was less than 2\%, however those results where discarded. Therefore, we can see that $\mathds{E}[\tilde{d}_e] = d_e, \forall e \in E$, except we control the variance by adjusting the min and max parameter.

\begin{table}[htb]
\begin{tabular}{c|ccc}
\textbf{MHE} & \textbf{V Nodes} & \textbf{Graph Nodes} & \textbf{\% Failure} \\ \hline
\textit{1}   & 0-3-9-5-11-13  & 6-E-A-C-D-6    & 0                   \\
\textit{2}   & 0-4-10-1-7-13  & 6-B-A-A-B-6    & \textbf{0.21675}    \\
\textit{3}   & 0-2-8-13       & 6-C-D-6        & \textbf{0.027}      \\
\textit{4}   & 0-6-12-13      & 6-E-C-6        & 0                   \\
\textit{5}   & 0-13           & 6-6            & 0
\end{tabular}
\caption{Deterministic results.}
\label{table1:det_res}
\end{table}

Table \ref{table1:det_res} presents the first set of results by solving the example using formulation \eqref{form:FS}. Second column shows nodes as described using the $V$ set, while the third one shows its correspondence with the true node in the graph. The final columns shows an out-of-sample analysis of failure: we created 1000 testing sets of random traveling times and we tested if the preset path was able to comply with the time windows. The total traveling distance between all MHE was 320 meters.

\begin{table}[htb]
\begin{tabular}{c|ccc}
\textbf{MHE} & \textbf{Tasks}     & \textbf{Nodes}  & \textbf{\% Failure} \\ \hline
\textit{1}   & 0-13               & 6-6             & 0                   \\
\textit{2}   & 0-6-12-13          & 6-E-C-6         & 0                   \\
\textit{3}   & 0-2-8-13           & 6-C-D-6         & \textbf{0.028}      \\
\textit{4}   & 0-3-9-13           & 6-E-A-6         & 0                   \\
\textit{5}   & 0-1-7-4-10-5-11-13 & 6-A-B-B-A-C-D-6 & \textbf{0.008667}
\end{tabular}
\caption{Stochastic results with $\alpha = 0.$}
\label{table2:sto_res}
\end{table}

Table \ref{table2:sto_res} presents the same analysis using formulation \eqref{form:SPM} and $\alpha = 0$. We can see a few differences in the paths chosen for each MHE, however the interest focuses on the small differences in terms of the pathing selected for each MHE. Moreover, thanks to the small changes we can check that the probability of failure in the out-of-sample analysis was never greater than 2\%, in contrast with the deterministic solution where its maximum was over 20\%. The usage of $\alpha$ set to 0 was to study the robust solution obtained, while solving the problem in a fraction of the time, since there is no need of creating more constraints than those for the supremum of the random traveling times sampled. Any value larger than 0 presented the same results, albeit larger solving times due to the integer complexity. The total traveling distance between all MHE was also 320 meters.

\section{Conclusion}

Production logistics systems exist in all production systems. And production logistics connects the components of the production system. Therefore, in order to operate the production system efficiently, the production logistics system must be operated efficiently, and in order to optimize the production system, the production logistics system must also be optimized. In this paper, we showed that the production system and the production logistics system are connected and proposed a method for optimizing the production logistics plan using stochastic variables.

Optimization models always give the best answer, and unlike heuristics, they guarantee an optimal solution. However, it is also true that defining the problem and constructing the model is time-consuming and needs effort. The method proposed in this paper does not completely solve the production logistics planning problem. Because it didn't organize the problem by integrating the sequencing of logistics tasks and assigning tasks to MHEs. But still, our proposed method is useful for drawing the big picture. Methods such as machine learning are more useful for dynamic environments where the problem is small and frequent updates are required. Nevertheless, the proposed method in our paper provides schedulers with a new perspective and a solution based on a mathematical model rather than an intuitive solution. In addition, factors such as sustainability can be included in the objective function, rather than simply considering efficiency.

From the early results, we see great opportunity in pursuing further research. The robust model was able to obtain a large reduction in the failure probability while maintaining the same \emph{costs} (i.e. total traveling distance) of its deterministic counterpart. We believe a more complex study case will present further difficulty and challenges, albeit more interesting results and analysis.

\bibliography{biblio}

\begin{thebibliography}{10}

\bibitem{mourtzis2020simulation}
Dimitris Mourtzis.
\newblock Simulation in the design and operation of manufacturing systems:
  state of the art and new trends.
\newblock {\em International Journal of Production Research}, 58(7):1927--1949,
  2020.

\bibitem{alexopoulos2018industrial}
Kosmas Alexopoulos, Konstantinos Sipsas, Evangelos Xanthakis, Sotiris Makris,
  and Dimitris Mourtzis.
\newblock An industrial internet of things based platform for context-aware
  information services in manufacturing.
\newblock {\em International Journal of Computer Integrated Manufacturing},
  31(11):1111--1123, 2018.

\bibitem{huang2017collaborative}
Houfeng Huang, Qing Ling, Wei Shi, and Jinlin Wang.
\newblock Collaborative resource allocation over a hybrid cloud center and edge
  server network.
\newblock {\em Journal of Computational Mathematics}, 35(4), 2017.

\bibitem{gilboa2006distributed}
Arnon Gilboa, Amnon Meisels, and Ariel Felner.
\newblock Distributed navigation in an unknown physical environment.
\newblock In {\em Proceedings of the fifth international joint conference on
  Autonomous agents and multiagent systems}, pages 553--560, 2006.

\bibitem{bhattacharya2011distributed}
Subhrajit Bhattacharya and Vijay Kumar.
\newblock Distributed optimization with pairwise constraints and its
  application to multi-robot path planning.
\newblock {\em Robotics: Science and Systems VI}, 177, 2011.

\bibitem{chouhan2015dmapp}
Satyendra~Singh Chouhan and Rajdeep Niyogi.
\newblock Dmapp: A distributed multi-agent path planning algorithm.
\newblock In {\em Australasian Joint Conference on Artificial Intelligence},
  pages 123--135. Springer, 2015.

\bibitem{banyai2019smart}
{\'A}gota B{\'a}nyai, B{\'e}la Ill{\'e}s, Elke Glistau, Norge Isaias~Coello
  Machado, P{\'e}ter Tam{\'a}s, Faiza Manzoor, and Tam{\'a}s B{\'a}nyai.
\newblock Smart cyber-physical manufacturing: Extended and real-time
  optimization of logistics resources in matrix production.
\newblock {\em Applied Sciences}, 9(7):1287, 2019.

\bibitem{liao2019novel}
Wenzhu Liao and Tong Wang.
\newblock A novel collaborative optimization model for job shop
  production--delivery considering time window and carbon emission.
\newblock {\em Sustainability}, 11(10):2781, 2019.

\bibitem{li2018new}
Guomin Li, Bing Zeng, Wei Liao, Xinyu Li, and Liang Gao.
\newblock A new agv scheduling algorithm based on harmony search for material
  transfer in a real-world manufacturing system.
\newblock {\em Advances in Mechanical Engineering}, 10(3):1687814018765560,
  2018.

\bibitem{yao2020improving}
Fengjia Yao, Bugra Alkan, Bilal Ahmad, and Robert Harrison.
\newblock Improving just-in-time delivery performance of iot-enabled flexible
  manufacturing systems with agv based material transportation.
\newblock {\em Sensors}, 20(21):6333, 2020.

\bibitem{otto2017reducing}
Alena Otto and Olga Batta{\"\i}a.
\newblock Reducing physical ergonomic risks at assembly lines by line balancing
  and job rotation: A survey.
\newblock {\em Computers \& Industrial Engineering}, 111:467--480, 2017.

\bibitem{lodree2009taxonomy}
Emmett~J Lodree~Jr, Christopher~D Geiger, and Xiaochun Jiang.
\newblock Taxonomy for integrating scheduling theory and human factors: Review
  and research opportunities.
\newblock {\em International Journal of Industrial Ergonomics}, 39(1):39--51,
  2009.

\bibitem{bechtsis2017sustainable}
Dimitrios Bechtsis, Naoum Tsolakis, Dimitrios Vlachos, and Eleftherios Iakovou.
\newblock Sustainable supply chain management in the digitalisation era: The
  impact of automated guided vehicles.
\newblock {\em Journal of Cleaner Production}, 142:3970--3984, 2017.

\bibitem{halawa2020introduction}
Farouq Halawa, Husam Dauod, In~Gyu Lee, Yinglei Li, Sang~Won Yoon, and
  Sung~Hoon Chung.
\newblock Introduction of a real time location system to enhance the warehouse
  safety and operational efficiency.
\newblock {\em International Journal of Production Economics}, 224:107541,
  2020.

\bibitem{battarra2014chapter}
Maria Battarra, Jean-Fran{\c{c}}ois Cordeau, and Manuel Iori.
\newblock Chapter 6: pickup-and-delivery problems for goods transportation.
\newblock In {\em Vehicle Routing: Problems, Methods, and Applications, Second
  Edition}, pages 161--191. SIAM, 2014.

\bibitem{kleywegt2002sample}
Anton~J Kleywegt, Alexander Shapiro, and Tito Homem-de Mello.
\newblock The sample average approximation method for stochastic discrete
  optimization.
\newblock {\em SIAM Journal on Optimization}, 12(2):479--502, 2002.

\bibitem{luedtke2008sample}
James Luedtke and Shabbir Ahmed.
\newblock A sample approximation approach for optimization with probabilistic
  constraints.
\newblock {\em SIAM Journal on Optimization}, 19(2):674--699, 2008.

\bibitem{pagnoncelli2009sample}
Bernardo~K Pagnoncelli, Shabbir Ahmed, and Alexander Shapiro.
\newblock Sample average approximation method for chance constrained
  programming: theory and applications.
\newblock {\em Journal of optimization theory and applications},
  142(2):399--416, 2009.

\end{thebibliography}
\bibliographystyle{unsrt}
\end{document}